%% file: calculs0.tex
\documentclass{amsart}
\usepackage[title,titletoc,toc]{appendix}
\usepackage{ifthen,amsthm,srcltx,amsopn,amssymb,amsfonts}
\usepackage{amssymb,color}
\usepackage{colordvi,fancyhdr}
\usepackage{pdfpages}
\definecolor{Gris}{cmyk}{0.1,0.1,0.1,.75}
\usepackage[colorlinks=true,citecolor=Gris,linkcolor=Gris,filecolor=Gris,urlcolor=Gris]{hyperref}
\usepackage{asymptote}
\usepackage{amscd}
\usepackage[english]{babel}
\RequirePackage{graphics,epsfig,psfrag}
\usepackage[T1]{fontenc}
\usepackage[utf8x]{inputenc}
\usepackage{delarray,graphicx}
\usepackage{comment}
\usepackage{tikz}
\usepackage[all,cmtip]{xy}
\usepackage{graphicx}

\usepackage{amssymb,latexsym}

\usepackage{amsthm}
\usepackage{eucal}

\usepackage{multirow}
\usepackage{lscape}
\usepackage{array}

\newcommand{\showcomments}{no}
\renewcommand{\showcomments}{no}

\newsavebox{\commentbox}
{\ifthenelse{\equal{\showcomments}{no}}%
{\footnotemark
        \begin{lrbox}{\commentbox}
        \begin{minipage}[t]{1.25in}\raggedright\sffamily\tiny
        \footnotemark[\arabic{footnote}]}
{\begin{lrbox}{\commentbox}}}%
{\ifthenelse{\equal{\showcomments}{no}}%
{\end{minipage}\end{lrbox}\marginpar{\usebox{\commentbox}}}
{\end{lrbox}}}

\title[]{Representations of fundamental groups of 3-manifolds into $\PGL(3,\bC)$: Exact computations in low complexity. }

\author{E. Falbel, P.-V. Koseleff and F. Rouillier}
\thanks{This work was supported in part by the ANR through the project "Structures G\'eom\'etriques et Triangulations".}
\address{Institut de Math\'ematiques de Jussieu \\
CNRS UMR 7586 and INRIA EPI-OURAGAN \\
Universit\'e Pierre et Marie Curie \\
4, place Jussieu 75252 Paris Cedex 05, France \\}
\email{}
\urladdr{}

\DeclareFontFamily{OT1}{rsfs}{}

\DeclareFontShape{OT1}{rsfs}{n}{it}{<-> rsfs10}{}
\DeclareMathAlphabet{\mathscr}{OT1}{rsfs}{n}{it}

\DeclareFontFamily{OT1}{rsfs}{}

\DeclareFontShape{OT1}{rsfs}{n}{it}{<-> rsfs10}{}
\DeclareMathAlphabet{\mathscr}{OT1}{rsfs}{n}{it}
\newcommand{\Q}{\mathbb{Q}}

\newcommand{\R}{\mathbb{R}}

\swapnumbers

\newtheorem{proposition}[subsection]{Proposition}

\newtheorem*{theorem*}{Theorem}
\newtheorem*{appl*}{Application}

\theoremstyle{definition}

\theoremstyle{remark}

\newcommand\bZ{\mathbb{Z}}
\newcommand\bC{\mathbb{C}}
\newcommand\bR{\mathbb{R}}

\newcommand\SL{\textrm{SL}}
\newcommand\PGL{\textrm{PGL}}
\newcommand\PU{\textrm{PU}}
\newcommand\Fl{\mathcal{F}l}
\renewcommand\P{\mathbb{P}}

\numberwithin{equation}{subsection}

\newcommand{\RP}{\mathbb{RP}}
\newcommand{\CP}{\mathbb{CP}}
\newcommand{\cal}{\mathcal}

\newcommand{\PSL}{\mathrm{PSL}}
\def\frac#1#2{{\textstyle{{#1} \overwithdelims.. {#2}}}}
\def\Frac#1#2{{\displaystyle{{#1} \overwithdelims.. {#2}}}}
\def\cZ{{\cal Z}}
\def\cS{{\cal S}}
\def\Mod#1{{\, (\mathrm{mod}\, #1)}}
\def\abs#1{\left \vert #1 \right \vert}
\def\t#1{\hbox{}^{t}\! #1}
\def\Re{\mathrm{Re}\,}
\def\Im{\mathrm{Im}\,}
\def\adots{\mathinner{\mkern2mu\raise1pt\hbox{.}
\mkern3mu\raise4pt\hbox{.}\mkern1mu\raise7pt\hbox{.}}}

\begin{document}

\begin{abstract}
In this paper we are interested in computing representations
of the fundamental group of a 3-manifold  into $\PGL(3,\bC)$ (in particular in $\PGL(2,\bC),
\PGL(3,\bR)$ and $\PU(2,1)$).  The representations are obtained
by gluing decorated tetrahedra of flags as in \cite{falbeleight, BFG}.  We list complete computations (giving 0-dimensional or 1-dimensional solution sets (for unipotent boundary holonomy)
for the first complete hyperbolic non-compact manifolds with finite volume which are obtained gluing less than three tetrahedra with a description of the computer methods used to find them.
The methods we use work for non-unipotent boundary holonomy as shown in some examples.
\end{abstract}

\maketitle
\begin{center}
\parbox{12cm}{\small
\tableofcontents
}
\end{center}
\input calculs0-I.tex

\input calculs0-II.tex
\input calculs0-III.tex

\bibliography{bibli-1}
\bibliographystyle{plain}

\end{document}

%% file: calculs0-I.tex
\section{Introduction}

In this paper we are interested in obtaining representations
of the fundamental group of a 3-manifold  into $\PGL(3,\bC)$.  They will be defined via certain topological triangulations with additional geometric data carried by the
0-skeleton of the triangulation.

The most important example being hyperbolic geometry we will mainly consider 3-manifolds which carry a complete hyperbolic structure.  That is, complete riemannian manifolds with  constant negative curvature equal to $-1$.

 In order to simplify the description we first treat open manifolds
which are the interior of manifolds with boundary.  If the hyperbolic manifold is complete and of finite volume one shows that its ideal boundary is a union of tori.  In that case we consider ideal triangulations such that the 0-skeleton is in the ideal boundary.  That is, identifying each boundary to a point, we consider a triangulation of the quotient space such that the 0-skeleton
coincides with the set of boundary points.

In hyperbolic space the ideal tetrahedra are described by a configuration of four points in the boundary of hyperbolic space, that is, $\CP^1$.  The cross-ratio is sufficient
to parameterize ideal tetrahedra up to the isometry group.  Ideal triangulations were used by Thurston to obtain examples of hyperbolic manifolds (see \cite{Thurston}).  A complete hyperbolic structure
is obtained on the manifold once a system of equations on the cross-ratio coordinates of those ideal tetrahedra are solved.  A very successful computational tool, SnapPea, was developed by Jeff Weeks (see its very readable description in \cite{Weeks}) which solves numerically the system of equations.

Other geometric structures in 3-manifolds are associated to subgroups of $\PGL(3,\bC)$, namely, flag structures which are associated to $\PGL(3,\bR)$ and spherical Cauchy-Riemann structures (CR structures)
which are associated to the group $\PU(2,1)$. These geometric structures are not well understood, in particular one does not know which 3-manifolds can carry one of them.

The first step in trying to find a structure is to obtain a representation of the fundamental group in one of these groups (see \cite{falbeleight} for the case $\PU(2,1)$).   We will use the raw data of triangulation of 3-manifolds as obtained in SnapPea.  Then we proceed as in \cite{BFG}; to a triangulation we associate a decoration
given by a choice of flags at each vertex of the tetrahedra.  We set equations imposing that the gluing of the tetrahedra is compatible with the decoration of flags.  The solution of the equations would give immediately a representation of the fundamental group of the 3-manifold into $\PGL(3,\bC)$ (called the holonomy representation).  Imposing that the holonomy restricted to the boundary tori be unipotent is a natural condition which is analogous to the condition of completeness in the hyperbolic case.  But more general conditions are very important too.  In particular, a condition, which is appropriate to the CR case is that the eigenvalues be of absolute value one.
Indeed, parabolic boundary conditions in the CR case imply that condition and it would be interesting to determine if there exist CR structures corresponding to  boundary holonomies given by parabolic
 conditions. Our methods to solve these equations are the same and we give in  section \ref{sec:m004ni}  solutions to some  non-unipotent
systems.  More general boundary conditions could be treated by our methods.

A generalization of the gluing equations for higher dimensions was described in \cite{GGZ}, and a simplified set of equations describing the particular case of representations into $\SL(n,\bC)$ is also described in \cite{GTZ1}.
These equations are based on decorations of tetrahedra by affine flags (see also the $a$-coordinates for affine flags  in \cite{BFG}) but, as there were initially designed for finding representations into $\SL(n,\bC)$, they miss solutions related to  $\PGL(n,\bC)$ representations obtained via projective flags. However, the introduction of a cocycle and its computation gives access to all boundary unipotent representation in $\PGL(n,\bC)$.  This interesting method has the advantage to
introduce quadratic equations as systems to be solved but is not adapted to non unipotent decorations and thus to describe other representations.
On the other hand, it is not clear whether the computations are actually simpler in the affine flag case as shown by limitations in the computations in \cite{Ptolemy} (even for $\PGL(3,\bC)$).  Indeed, the use of Gröbner bases discards immediately the initial structure of the equations (symetries, sparsity, etc.) and the small degrees of the initial equations does not prevent the appearance of very large degrees during the computations. Even if some classes of examples are known to have a good behaviour (see \cite{HL2011} for $0$-dimensional systems), up to our knowledge, there does not exist any criteria on the initial equations that might help to decide if such algorithm will be well behaved or not.

In this paper we deal with methods to solve the equations and obtain a list of solutions for manifolds with low complexity.  In particular we obtain all
solutions for ideal triangulations of complete cusped hyperbolic manifolds with less than four tetrahedra.

In the first sections, in order to make the paper self-contained,  we review results in \cite{BFG} (see also \cite{falbeleight} for the $\PU(2,1)$ case).  That contains the parametrization of decorated tetrahedra (configurations of four flags in $\CP^2$),  the  description of decorated ideal triangulations, the compatibility equations
which will lead us to a system of equations and the computation of the holonomy representation.  The special case we deal mostly in this paper has unipotent boundary holonomy.

In section \ref{tools} we describe the methods used to solve the system of equations.

In the remaining sections we describe several important examples which illustrate the methods and the results.
 The complete census up to three tetrahedra is shown in Table \ref{DesSols} (more details can be obtained in our webpage \href{https://who.rocq.inria.fr/Fabrice.Rouillier/SGT/Home_page.html}{.../SGT}).

A very important observation which came out from the examples is that solutions of the gluing equations with unipotent boundary holonomy might not be 0-dimensional
even in low complexity.  Computations for representations in $\PGL(2,\bC)$ show that, for cusped hyperbolic manifolds obtained with  four tetrahedra, only two have a one parameter family component (there are no examples with less than four tetrahedra), but we do not have an efficient method to decide when the solution will have
a positive dimensional component.
We verified that the one dimensional components (in all examples) have at most a finite number of $\PSL(2,\bC)$ solutions. Moreover, they are all in $\PSL(2,\bR)$ so the volumes presented in the table concern indeed all $\PSL(2,\bC)$ representations arising from the given triangulation.
On the other hand, $\PU(2,1)$ solutions might appear in families (of real dimension one or two) inside 1-dimensional  components (see the worked example m003 in \cite{BFGKR} where there is a whole complex 1-dimensional component consisting of CR solutions).   In any case, all these families arise from degenerate representations.  Indeed, either
they give rise to reducible representations or they have trivial boundary holonomy.

In the case of the figure eight knot the non-hyperbolic representations in $\PGL(3,\bC)$ were  obtained in \cite{falbeleight}.  They happen to be discrete representations in $\textrm{\PU}(2,1)$ and two of them give rise to spherical structures on the complement of the figure eight knot (see \cite{DF,FW}).  In the case of the Whitehead link complement, R. Schwartz (\cite{Schwartz}) was able to analyse one  representation (obtained in a completely independent way as a subgroup of $\textrm{\PU}(2,1)$ generated by reflections) and showed that the complement of the link has a spherical CR structure.  Other representations of the fundamental group of the complement of the Whitehead link obtained similarly using reflection groups were analysed in \cite{PW}.  In their case the representations form a one parameter family parameterized by the boundary holonomy. We obtain here independently the unipotent representation which is  the holonomy of a CR structure on the complement of the Whitehead link.

It turns out that for the examples of CR structures obtained until now, the boundary holonomy of the $\PU(2,1)$ representations are abelian of rank one.  This fact made us  try
to chase  representations into $\PU(2,1)$ with rank one boundary holonomy in the hope that they will correspond to CR structures.  We enumerate all these solutions.
  For the ideal triangulation of the Whitehead link complement with four tetrahedra there is only one unipotent
representation with rank one  boundary holonomy (it is the same one obtained independently in \cite{PW}).

Above all,  the computations encourage us to study CR structures.  There is a great proportion of $\PU(2,1)$ representations obtained among all $\PGL(3,\bC)$ representations
and some of them might correspond to holonomies of geometric structures.
On the other hand it is a challenging problem to decide if a 3-manifold has a geometric structure.  It is not clear when the representations obtained
here are discrete (for general discrete subgroups of $\PGL(3,\bC)$ see \cite{CNS}) .  On the other hand, a criterium for
 rigidity of representations is given in \cite{BFGKR} and a discussion of generic discreteness of  the boundary holonomy is given in \cite{G}.  In this paper we don't address these questions.

Another interesting observation from  Table \ref{DesSols} is the fact that the hyperbolic volume is the maximum among the volumes of all unipotent $\PGL(3,\bC)$ representations.  This was shown recently using methods of bounded cohomology in \cite{BBI}.   On the other hand, $\PU(2,1)$ representations have null volume
by cohomological reasons (see \cite{FWang}).

A final remark is the fact that no numerical method is available to solve these equations. It is a remarkable fact that, due to Mostow's rigidity, one knows that  there is  at most
one complete hyperbolic structure on a manifold.
This structure corresponds to the unique solution with positive imaginary parts of the variables associated to a triangulation.
The numerical scheme introduced in SnapPea should converge to this solution.
Then, the LLL-algorithm makes the method very efficient in finding exact solutions corresponding to the hyperbolic structures.
On the other hand, there is no knowledge of the number or location of solutions to the equations for representations in $\PGL(3,\bC)$. Therefore, there is no numerical scheme to obtain solutions and we are, for the moment, condemned to solve the equations exactly.

We thank M. Thistlethwaite for introducing us to SnapPea and making available his list of cusped 3-manifolds with particularly simple generators
for the boundary fundamental group.  We used his list for our computations.  We also thank N. Bergeron, M. Deraux, A. Guilloux,
 M. Thistlethwaite
and S. Tillmann for all the stimulating discussions.

\section{Three geometric structures on three manifolds}

Geometric structures on three manifolds have been extensively studied.  The usual setting is an action of
a Lie group $G$ on a homogeneous 3-manifold $X$.  An $(X,G)$ structure on a 3-manifold $M$
 being a family of charts $\phi_i : U_i\rightarrow X$
from open subsets forming a cover of $M$ with transition functions $g_{ij}$ (given by $\phi_j=g_{ji}\circ \phi_i $) in the Lie group $G$.

Hyperbolic structures, that is, $({\mathbb H}_{\bR}^{3},\PSL(2,\bC))$ structures (where ${\mathbb H}_{\bR}^{3}$ is hyperbolic 3-space) were shown by Thurston to be very important and his theory
made possible, as a far reaching consequence of his ideas, to solve Poincar\'e's conjecture.   The boundary of the 3-dimensional hyperbolic space can be identified to $\CP^1$.  By embedding $\CP^1$ into $\CP^2$ as a conic
we observe that a point in $\CP^1$ defines a point and a line containing it in $\CP^2$.  Namely the point obtained by the embedding and the tangent line to
the embedding passing through that point.   This justifies considering configuration of flags associated to triangulations (recall that a pair consisting of a point and a projective line containing the point is called a flag).  We will associate to each vertex of a tetrahedron
a flag in the following section.

 CR geometry is modelled on the sphere ${\mathbb S}^3$ equipped with a natural $\mathrm{PU}(2,1)$ action.  More precisely, consider the group $\mathrm{U}(2,1)$ preserving the Hermitian form
$\langle z,w \rangle = w^*Jz$ defined on ${\bC}^{3}$ by
the matrix
$$
J=\left ( \begin{array}{ccc}

                        0      &  0    &       1 \\

                        0       &  1    &       0\\

                        1       &  0    &       0

                \end{array} \right )
$$
and the following cones in ${\mathbb C}^{3}$;
$$
        V_0 = \left\{ z\in {\bC}^{3}-\{0\}\ \ :\ \
 \langle z,z\rangle = 0 \ \right\},
$$
$$
        V_-   = \left\{ z\in {\bC}^{3}\ \ :\ \ \langle z,z\rangle < 0
\ \right\}.
$$
Let $\pi :{\bC}^{3}\setminus\{ 0\} \rightarrow \mathbb{CP}^{2}$ be the
canonical projection.  Then
${\mathbb H}_{\bC}^{2} = \pi(V_-)$ is the complex hyperbolic space and its boundary is
$$
\partial{\mathbb H}_{\bC}^{2} ={\mathbb S}^3= \pi(V_0)=\{ [x,y,z]\in \mathbb{CP}^{2}\ |\ x\bar z+ |y|^2+z\bar x=0\ \}.
$$
The group of biholomorphic transformations of ${\mathbb H}_{\bC}^{2}$ is then
$\mathrm{PU}(2,1)$, the projectivization of $\mathrm{U}(2,1)$.  It acts on ${\mathbb S}^3$ by
CR transformations.
An element $x\in {\mathbb S}^3$ gives rise to a flag in ${\CP}^{2}$ where the line corresponds to the unique complex line tangent to ${\mathbb S}^3$ at $x$.

A third real 3-dimensional geometry is the geometry of real flags in $\bR^3$.  That is the geometry of the space of all couples $[p,l]$ where $p\in \RP^2$ and $l$ is a projective line
containing $p$.  The space of flags is identified to the quotient
$$
\SL(3,\bR)/B
$$
where $B$ is the Borel group of all upper triangular matrices.

In the next section we describe the space of flags in $\CP^2$ which will be the common framework to describe the three geometries based on
$\PSL(2,\bC), \PU(2,1)$ and $\PGL(3,\bR)$.

\section{Flag tetrahedra}\label{ss:coord-para}

In this section we  recall the parametrization  of configurations of four flags in the projective space $ \CP^2$ (more details can be seen in \cite{BFG}).
Let $V=\bC^3$. A flag in $V$ is usually seen as a line and a plane, the line belonging to the plane.  Using the dual vector space $V^*$ and the projective spaces $\P(V)$ and $\P(V^*)$,
define the spaces of  {\it flags} $\Fl$ by the following:
\begin{eqnarray*}
 \Fl & = & \{([x],[f])\in \P(V)\times \P(V^*) \, | \, f(x)=0\}.
\end{eqnarray*}

The natural action of $\SL(3,\bC)$ on $\Fl$
makes us identify the space of flags  with the homogeneous space $\SL(3,\bC)/ B$, where $B$ is the Borel subgroup of upper-triangular matrices in $\SL(3,\bC)$.
A {\it generic configuration of  flags} $([x_i],[f_i])$, $1\leq i\leq n+1$
is given by $n+1$ points $[x_i]$ in general position and $n+1$ lines $f_i$
in $\P(V)$ such that $f_j(x_i)\neq 0$ if $i\neq j$.
A configuration of ordered points in $\P(V)$ is said to be in {\it general position} when they are all distinct and no three points are contained in the same line.


Up to the action of $\PGL(3,\bC)$, a generic configuration of three flags $([x_i],[f_i])_{1\leq i \leq 3}$ has  only one invariant given by  the triple ratio
$$
X=\Frac{f_1(x_2)f_2(x_3)f_3(x_1)}{f_1(x_3)f_2(x_1)f_3(x_2)}\in \bC^{\times}
$$

\subsection{Coordinates for a tetrahedron of flags}

\begin{figure}[!ht]
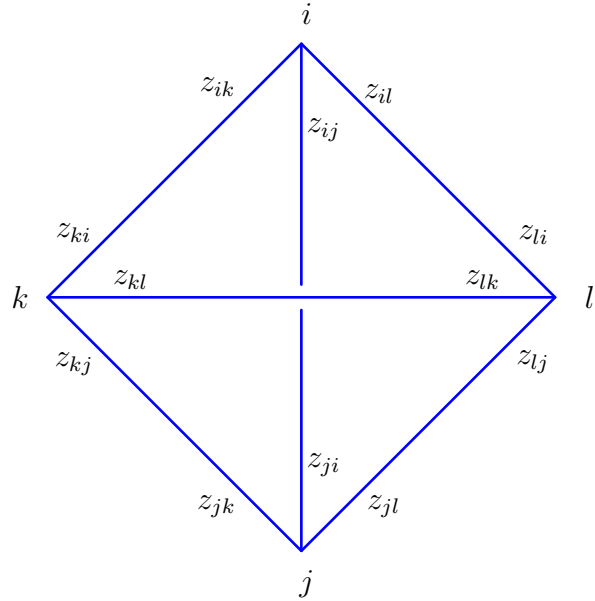
 \label{tetra}
\begin{center}
\begin{asy}
  // 
  // 
  size(8cm);
  defaultpen(1);
  usepackage("amssymb");
  import geometry;

  // 
  point o = (8,0);
  pair ph = (8,0);
  pair pv = (4,4);

  // 
  draw(o -- o+pv -- o+ph -- o+ph-pv --cycle, blue);
  draw(o +pv --o+ ph/2+(0,.2), blue);
  draw(o+ph-pv -- (12,-.2), blue);
  draw(o -- o+ph , blue);

  // 
  // dot(o,black+5);
  // 
  label("{ $z_{kl}$}", o, 6*dir(15));
  label(" $ z_{kj}$", o, 6*dir(-75));
  label(" $z_{ki}$", o, 6*dir(75));
  label("{ $k$}", o, 2*W);

  label("{$z_{lk}$}", o+ph, 6*dir(180-15));
  label(" $ z_{lj}$", o+ph, 6*dir(-110));
  label(" $z_{li}$", o+ph, 6*dir(110));
  label("{ $l$}", o+ph, 2*E);

  label("{$z_{il}$}", o+pv, 8*dir(-35));
  label(" $ z_{ij}$", o+pv, 8*dir(-80));
  label(" $z_{ik}$", o+pv, 8*dir(-150));
  label("{ $i$}", o+pv, 2*N);

   label("{$z_{jl}$}", o+ph-pv, 8*dir(30));
  label(" $ z_{ji}$", o+ph-pv, 8*dir(80));
  label(" $z_{jk}$", o+ph-pv, 8*dir(150));
  label("{ $j$}", o+ph -pv, 2*S);

   // 
  // pair p1 = (4, 2);
  // pair p2 = (7, -1);
  // guide H = p1{E}..{S}p2;
  // draw(H,red);

\end{asy}
\caption{The cross-ratio coordinates (or $z$-coordinates).} 
\end{center}
\end{figure}

We recall the parametrization used in \cite{BFG}.  We refer to Figure \ref{tetra} which displays the coordinates.

Let $([x_i],[f_i])_{1\leq i \leq 4}$ be a generic tetrahedron.
 We define a set of 12 coordinates on the edges of the tetrahedron (one for each oriented edge).
To define the coordinate $z_{ij}$ associated to the edge $ij$, we first define $k$ and $l$ such that the permutation $(1,2,3,4) \mapsto (i,j,k,l)$ is even. The pencil of (projective) lines through the point $x_i$ is a projective line $\P_1(\bC)$. We have four points in this projective line: the line $\textrm{ker}(f_i)$ and the three lines through $x_i$ and one of the $x_l$ for $l\neq i$. We define $z_{ij}$ as the cross-ratio of four flags
by
$$z_{ij} := [\textrm{ker}(f_i),(x_ix_j),(x_ix_k),(x_ix_l)].$$
Note that we follow the convention that the cross-ratio of four points
$x_1, x_2 , x_3 , x_4$ on a line is the value at $x_4$ of a projective coordinate taking
value $\infty$ at $x_1$, $0$ at $x_2$, and $1$ at $x_3$. So we employ the formula
$$[x_1,x_2,x_3,x_4] :=\Frac{(x_1 -x_3)(x_2-x_4)}{(x_1 -x_4)(x_2-x_3)}$$ for the cross-ratio.

At each face $(ijk)$ (oriented as the boundary of the tetrahedron $(1234)$), the $3$-ratio is the opposite of the product of all cross-ratios ``leaving'' this face:
$$z_{ijk} := \Frac{f_i(x_j)f_j(x_k)f_k(x_i)}{f_i(x_k)f_j(x_i)f_k(x_j)}=- z_{il}z_{jl}z_{kl}.$$

This follows from a direct computation (see \cite{BFG}).
Observe that if the same face $(ikj)$ (with opposite orientation)  is common to a second tetrahedron $T'$ then $$z_{ikj}(T) = \Frac{1}{z_{ijk}(T')}.$$ 

Of course there are relations between the whole set of coordinates, namely,
the three cross-ratio leaving a vertex are algebraically related:
\begin{equation}
\begin{split}
& z_{ik} =  \Frac{1}{1-z_{ij}},\\
& z_{il} =   1-\Frac{1}{z_{ij}}.
\end{split}
\end{equation} \label{Cros}

Observe that the relations follow a cyclic order around each vertex which is defined by the orientation of the tetrahedron.
We also have
\begin{equation} \label{3-1}
z_{ij} \, z_{ik} \, z_{il} = -1.
\end{equation}
The next proposition shows that a tetrahedron is uniquely determined, up to the action of $\PGL(3,\bC)$, by four numbers.
One can pick a variable at each vertex.  In fact, the space of flags  is a complex manifold ($\PGL(3,\bC)/B$) and, therefore, configurations of flags have
a natural complex structure.
\begin{proposition}\cite{BFG}\label{parameters}
The space of generic tetrahedra is biholomorphic to $(\bC\setminus \{0, \, 1\})^4$.
\end{proposition}

One can use one  cross-ratio coordinate at each vertex, for instance $(z_{12},z_{21},z_{34},z_{43})$ or  $(z_{13},z_{21},z_{32},z_{43})$.

\subsection{Hyperbolic ideal tetrahedra}

An ideal hyperbolic tetrahedron is given by $4$ points on the boundary of $\mathbb {\mathbb H}_{\bR}^{3}$,
i.e. $\P_1(\bC)$. Up to the action of $\PSL(2,\bC)$, these points are in homogeneous coordinates
 $[0,1]$, $[1,0]$, $[1,1]$ and $[1,z]$ -- the complex number $z$ being the cross-ratio of these four points.

 We may embed $\P_1(\bC)$ into  $\P_2(\bC)$
via the Veronese map $h$, given in homogeneous coordinates by
\begin{eqnarray*}
 h : \left[x,y\right] & \mapsto &  [x^2,xy,y^2]
\end{eqnarray*}
and therefore this map and its derivative define a map  from $\P_1(\bC)$ to the variety of flags $\Fl$.

Let $T$ be the tetrahedron $h([0,1])$, $h([1,0])$, $h([1,1])$ and $h([1,z])$.  Its image in the
variety of flags given by the above map has coordinates (see \cite{BFG})
\begin{eqnarray}
 z_{12}(T) = z_{21}(T)=z_{34}(T)=z_{43}(T)=z.\label{eq:hyp}
 \end{eqnarray}
Conversely, given parameters satisfying the equations above, they define a unique hyperbolic tetrahedron.

\subsection{The CR case.}\label{CR}

Recall that an element $x\in {\mathbb S}^3$ gives rise to an element  $([x],[f])\in \Fl(\bC)$ where $[f]$ corresponds to the unique complex line tangent to ${\mathbb S}^3$ at $x$.
The following proposition  describes the space of
generic configurations of four points in ${\mathbb S}^3$.

\begin{proposition}[\cite{falbeleight,BFG}]\label{prop:generic}
Generic configurations (up to translations by $\mathrm{PU}(2,1)$) of four points
in ${\mathbb S}^3$ not contained in an $\bR$-circle are parametrized
by generic configurations of four flags
with coordinates $z_{ij}$,  $1\leq i\neq j\leq 4$ satisfying  the three complex equations
\begin{eqnarray}\label{eq:cr}
z_{ij}z_{ji}=\overline {z_{kl}}\, \overline{z_{lk}}
\end{eqnarray}
not all of them being real
and such that  $z_{ji}z_{ki}z_{li}\neq 1$ for each face with vertices $\{ j,k,l\}$.
\end{proposition}
The conditions (\ref{eq:cr}) together with relation (\ref{3-1}) imply that the triple ratio of each face satisfies
$
\vert z_{ijk} \vert = 1.
$
For the tetrahedron to be
CR we need to verify the condition   $z_{ijk}=-z_{ji}z_{ki}z_{li}\neq -1$ for each face with vertices $\{ j,k,l\}$
(the triple ratio is never $-1$ for a CR generic triple of flags).
\medskip\par
Conditions $z_{ij}z_{ji}=\overline {z_{kl}}\, \overline{z_{lk}} \in \bR$ might describe other
configurations of four flags which are not CR. From $z_{ijk}z_{ilj} = \Frac{z_{lk}z_{kl}}{z_{ij}z_{ji}}$,
we deduce that in this case
$z_{ijk}z_{ilj}=1$ and therefore $z_{ijk} = z_{ilj} = z_{ikl} = z_{jkl}=\pm 1$.  If the cross-ratio invariants are all real, either
one obtains configurations contained in an $\bR$-circle (they coincide with the real hyperbolic ones with real cross ratios, that is,
$ z_{12}(T) = z_{21}(T)=z_{34}(T)=z_{43}(T)=x\in \bR- \{0,1\}$)
or one of the triple ratios is $-1$ and, in that case, they are not CR.


\subsection{The $\SL(3,\bR)$ case}\label{SLR}

Clearly, a configuration contained in the space of real flags is characterized by having all its
invariants $z_{ij}$ real.  Observe that the three cases intersect precisely for configurations corresponding to degenerate real
hyperbolic configurations.

\section{Ideal triangulations by flag tetrahedra}

We use the definition of an ideal triangulation of a 3-manifold as it is used usually in computations with SnapPea.
It is a union of 3-simplices $K= \bigcup_\nu T_\nu$ with face identifications (which are simplicial maps).  We let  $K^{(0)}$ be the union of vertices in
$K$.  The manifold $K- K^{(0)} $ is, topologically, the interior of a compact
manifold when the vertices $K^{(0)} $  are deleted. In fact, $K-K^{(0)}$
is an (open) $3$-manifold  that retracts onto a compact $3$-manifold with boundary $M$ (the boundary being the link of $K^{(0)}$).

Given such a compact oriented $3$-manifold $M$ with boundary, we call a triangulation $K$  as above an {\it ideal triangulation} of $M$.
A {\it parabolic decoration} of an ideal  triangulation is the data of a flag for each vertex (equivalently a map from the $0$-skeleton of the complex to $\Fl$).
A parabolic decoration
together with an ordering of the vertices of each $3$-simplex  equip each tetrahedron with a set of coordinates as defined in section \ref{ss:coord-para}.

From now on we fix $K$ a decorated oriented ideal triangulation of a $3$-manifold together with an ordering of the vertices of each $3$-simplex of $K$. Denote by $T_{\nu}$, $\nu =1 , \ldots , N$,  the tetrahedra of $K$ and
$z_{ij}(T_{\nu})$ the corresponding $z$-coordinates (or cross-ratio coordinates).  We impose the following compatibility conditions which will imply the existence of a well defined representation
of the fundamental group of the manifold $M$ into $\PGL(3,\bC)$.

\subsection{Consistency relations}
\label{sec-consistency}(cf. \cite{falbel,BFG})

\medskip

\emph{(Face equations)}
Let $T$ and $T'$ be two tetrahedra  of $K$ with a common face $(ijk)$ (oriented as a boundary of $T$), then $z_{ijk}(T)z_{ikj}(T')=1$.

\medskip

For a fixed edge   $e\in K$ let $T_{\nu_1}, \ldots , T_{\nu_{n_e}}$ be the $n_e$ tetrahedra in $K$ which contain an edge which projects unto  $e\in K$ (counted with multiplicity).  For each tetrahedron in $K$ as above we consider its edge $ij$ corresponding to $e$.

\medskip

\emph{(Edge equations)} $z_{ij}(T_{\nu_1}) \cdots z_{ij}(T_{\nu_{n_e}})=z_{ji}(T_{\nu_1}) \cdots z_{ji}(T_{\nu_{n_e}})=1$.

\medskip

The face equations are clearly necessary in order to match a triple of flags from one face to a triple of another face.  The edge equations follow by considering the 1-dimensional projective space of complex lines at each vertex.  Indeed, take a vertex whose associated flag is $(p_1,l_1)$ on the edge $[(p_1,l_1),(p_2,l_2)]$.  Consider the projective space of all lines through $p_1$.  The line $l_1$ will be identified to $\infty$ and the line $[p_1,p_2]$ to 0.  All the flags (not coinciding to $\infty$ or 0) in the tetrahedra
having this edge in common give rise to an ordered sequence of points in this 1-dimensional projective space.   Imposing that the the point corresponding to the last vertex of the last tetrahedron coincides with the point corresponding to the first vertex of the first tetrahedron amounts
precisely to the first edge condition.  Analogously, the second condition follows if we consider the projective space of lines at the other vertex of the edge.

One should be aware that the compatibility conditions do not imply immediately the existence of a geometric structure.
One should think of these as sufficient conditions for the existence of a 0-skeleton compatible with the side pairings.
Certainly it implies the existence of a representation of the fundamental group but one has yet to construct a compatible extension of
the side pairings to 3 simplices.   In the hyperbolic case one has the advantage of the existence of convex ideal
tetrahedra but in $\PU(2,1)$ there is no such canonical construction (see the discussion in \cite{falbeleight}).

\subsection{Volume}
\label{sec-volume}
We recall the definition of volume of a tetrahedron of flags in \cite{BFG}.
The {\it Bloch-Wigner dilogarithm} function is
\begin{align*}
D(x) & = \arg{(1-x)}\log{|x|}-\mathrm{Im} (\int_{0}^{x}\log{(1-t)}\frac{dt}{t}), \\
& = \arg{(1-x)}\log{|x|} + \mathrm{Im}( \ln_2 (x)).
\end{align*}
Here  $ \ln_2 (x)=-\int_{0}^{x}\log{(1-t)}\frac{dt}{t}$ is the dilogarithm function.
The function $D$ is well-defined and real analytic on $\bC-\{0,1\}$ and
extends to a continuous function on $\bC P^1$ by
defining $D(0) = D(1) = D(\infty) = 0$.

Given a tetrahedron of flags with coordinates  $\left(z_{1}(T), z_{2}(T),z_{3}(T),z_{4}(T)\right)$, we define its volume as
$$
{\rm{Vol}}(T)=\frac{1}{4}(D(z_1)+D(z_2)+D(z_3)+D(z_4)).
$$
The volume of a triangulation is the sum of the volumes of its tetrahedra.  If the tetrahedron is hyperbolic this definition coincides with
the volume of volume in hyperbolic geometry.  If the  tetrahedron is real its volume is zero.  Although CR tetrahedra have generically a non zero volume, if a triangulation is such that
all tetrahedra are CR  then
the total  volume is zero (\cite{FWang}).

\section{Holonomy of a decoration}
\label{sec-holonomy}
In this section we recall how to compute the holonomy of a decoration described above (see \cite{BFG} for more details).
A decoration of a triangulation of a manifold $M$ gives rise to a representation (called holonomy representation)
$$
 \rho : \pi_1(M,p_0)  \rightarrow \PGL(3,\bC).
$$
The base point is not important if we study representations up to conjugation.  In fact, each solution of the consistency equations gives
a conjugacy class of representations but a special choice of base point and side pairings  is needed  in explicit computations.
The idea is to follow a path in the fundamental group from face to face keeping track of changes using a coordinate system adapted to the faces.

%

In order to compute holonomies we first fix
a face (with ordered vertices) of one of the tetrahedra and a base point on  the face.  We then follow  paths
representing  the generators of the  fundamental group of the manifold which we decompose  into arcs contained in each tetrahedron
and which are transverse to their sides. For each arc we write the contribution to the holonomy
as is explained in the following.

Given a tetrahedron $(ijkl)$, the  arc going from face $(ijk)$ to face $(ijl)$ can be represented by a "left turn" arc
 around the edge $(ij)$.
We will also consider  permutations of a face as for example, the permutation
$(ijk)\rightarrow (jki)$.

To follow the change in the holonomy  along these arcs we associate to each face, that is, to a configuration of $3$ generic flags $([x_i ] , [f_i ] )_{1 \leq i \leq 3} $
with triple ratio $X$,  a projective
coordinate system   of $\CP^2$: take the one where the point $x_1=[1:0 : 0]^t$, $f_1=[0:0:1]$, $x_2=[0:0:1]^t$, $f_2=[1:0:0]$, the point $x_3$ has coordinates $[1:-1:1]^t$ and $f_3=[X :  X+1 : 1]$.

The cyclic permutation of the flags $(ijk)\rightarrow (jki)$ with triple ratio $X$ induces the coordinate change given
by the matrix
$$T(X) =  \left(
\begin{matrix}
X & X+1 & 1 \\
-X & -X & 0 \\
X & 0 &  0
\end{matrix} \right).$$

It is also useful to consider the transposition $ (ijk)\rightarrow (jik)$
which is given by the matrix
$$
I=\left(
\begin{matrix}
0& 0 & 1\\
0 & 1 & 0 \\
1 & 0 & 0
\end{matrix} \right).$$

Remark then that
$$
 T(\frac{1}{X} )=I\circ  T^{-1}(X)\circ I.
$$

If we have a tetrahedron of flags $(ijkl)$ with its $z$-coordinates, then the  basis
related to the triple $(ijl)$ is obtained in the  coordinate system related to the triple $(ijk)$
by the coordinate change given by the matrix
$$(ijk)\rightarrow (ijl) = \left(
\begin{matrix}
\Frac{1}{z_{ji}}& 0 & 0\\
0 & 1 & 0 \\
0 & 0 & z_{ij}
\end{matrix} \right).$$

\subsection{Representation of the fundamental group}

We  obtain the change of coordinate matrix of a path in the triangulation by decomposing it in the
two elementary steps described in the previous section and multiplying the change of coordinate matrices for each step from left to right.
Finally, to obtain the representation of the fundamental group we take the inverse of the coordinate change matrix for each path
representative of an element of the fundamental group.

One has to be careful if we want to obtain a representation in $\PU(2,1)$.
In general we need to conjugate the above representation by an element so that the group preserves a
chosen hermitian form.

Let 
$(i_0j_0k_0)$ (an ordered triple of points in a face of tetrahedron $T_0$)  be the base point and $Z_0 = z_{i_0j_0k_0}(T_0)$ be its 3-ratio.
We consider the transformation
$
A= \begin{pmatrix}
-\frac{1+Z_0}{Z_0} &  0 & 0 \\
0 & -1 & 0  \\
0 & 0 & 1 \\
\end{pmatrix}
$
such that
$
A \cdot
\left (\begin{matrix}\frac{1}{1+Z_0} \\1 \\1\end{matrix} \right)=
\left(\begin{matrix}1 \\-1\\1\end{matrix}\right).$
The hermitian form $J$ is preserved if we conjugate the group elements by $A$,
or equivalently the group preserves the hermitian form
$J_0 =\t\overline{A}^{-1} J A^{-1}$ with eigenvalues $1,  1/\abs{1+Z_0}, -1/\abs{1+Z_0}$.
We will show computations in a specific example in a future section.

\subsection{Boundary holonomy}
If we follow a path in a boundary torus of a triangulated manifold we will obtain, applying the procedure above, a coordinate change.  We obtain then a representation of the fundamental group of the boundary in $\PGL (3,\bC)$.
Each torus is associated to a vertex and therefore the representation of its fundamental
group fixes a flag.  Choosing the flag properly up to conjugation, we can arrange it so that the boundary representation is upper triangular.
   Given two generators of a torus group represented by paths $a$ and $b$ one can compute easily the eigenvalues of the holonomy
along them.

 We suppose a path is contained in the link of the vertex and is transverse to the edges of the triangulation induced by the tetrahedra.
Consider a path segment  at the link of vertex $i$ which is turning left around the edge $ij$, that is $(ijk)\to (ijl).$ Then, as before, the holonomy change of coordinate matrix is
$$
\left(
\begin{matrix}
\frac{1}{z_{ji}}& 0 & 0\\
0 & 1 & 0 \\
0 & 0 & z_{ij}
\end{matrix} \right).
$$
Consider the path segment at the link of vertex $i$ which is turning right around the edge $ij$, that is $$(ilj)\to (ikj).$$ Then the holonomy change of coordinate matrix is
$$
  \begin {pmatrix}
  \frac {z_{kl}z_{lk}}{z_{ij}}&
  -\frac{z_{kl}}{z_{ij}z_{lj}} - \frac{1}{z_{ki}}
&-
\frac {z_{il}}{z_{ki}}\\ \noalign{\medskip}0
&1&z_{il}\\ \noalign{\medskip}0&0&\frac{1}{z_{ij}}
\end {pmatrix}.
$$

 Both matrices being upper diagonal, to find the eigenvalues we need only to multiply the diagonal terms. The eigenvalues of the boundary holonomy are
 the diagonal elements obtained as above.  We will note $A^*$ and $A$ the two eigenvalues (the first and the third with the normalization as above so that the second eigenvalue be one) of the path $a$ and   $B^*$ and $B$ the two eigenvalues of the path $b$.

\subsection{Boundary unipotent holonomy equations}\label{sec:buh}
The boundary holonomy gives a representation of the fundamental group of the boundary tori into $\PGL(3,\bC)$.  If $t$ is the number of boundary tori, we note
$$
 \rho_b : (\bZ\oplus \bZ)^t \rightarrow \PGL(3,\bC).
$$
For knots we can define canonical generators of the boundary torus called meridian and longitude.  But, for simplicity, we will refer
to generators of boundary tori for other manifolds as well as meridians and longitudes, denoted $m$ and $l$.
The abelian group generated by the image of $\rho_b$ is called the peripheral group.
It turns out that an important information on  the representations is the rank of the peripheral group.  For instance, examples
of uniformizable CR manifolds were obtained precisely when the rank of the peripheral group is one (see \cite{Schwartz,DF,PW}).

We say the representation is $\sl unipotent$ if the boundary holonomy is unipotent.  That is, if the boundary holonomy representation can be given by matrices (up to scalar multiplication) of the form
$$ \left(
\begin{matrix}
1 & \star & \star\\
0 & 1 & \star \\
0 & 0 & 1
\end{matrix} \right ) .$$
That is, all eigenvalues are equal.

 In the following we will find unipotent representations associated to parabolic decorations.  That means that for each boundary torus
 we compute the eigenvalues  $A^*, A, B^*, B$ corresponding to generators $a$ and $b$ at each torus and impose the equations
 $$
 A^*= A= B^*= B=1.
 $$
 We remark that each eigenvalue is a rational function of the cross-ratio coordinates.

%% file: calculs0-II.tex
\section{Computational steps}\label{tools}

Our goal is to obtain the solutions of the system in the cross-ratio coordinates $z_{ij}$ and then
compute the holonomy representations.  We compute:
\begin{enumerate}
\item Solving the constraints on gluing decorated tetrahedra, which reduces to study a constructible set defined implicitly by
\begin{enumerate}
\item The consistency relations which define Edge and Face equations (\ref{sec-consistency}),
\item The holonomy equations (see section \ref{sec-holonomy}),
\item The cross-ratio equations (see \ref{Cros}),
\item The restrictions on the coordinate values, precisely $z_{i,j}$ must be different from $0$ and $1$.
\end{enumerate}
\item Distinguishing hyperbolic, real flags and CR solutions (see \ref{eq:cr})
\item Computing the resulting volumes (\ref{sec-volume})
\item Computing representations of the fundamental group (see section \ref{sec-holonomy}).
\end{enumerate}
Each system depends on the $12n$ variables
$\cZ = \{ z_{ij}(T)\}$, where $T$ denotes a tetrahedron in the triangulation (containing $n$ tetrahedra) and ${ij}$ one of its oriented edges.
\begin{itemize}
\item[---] $L_e$ (resp. $L_f$) denotes the set of the $2n$ edge (resp. face) equations;
\item[---] $L_h$ denotes the set of the polynomials defining the holonomy equations (their number is $4 t$, where $t$ is the number of tori in the boundary and in most cases treated here they correspond to unipotent conditions).
\item[---] $L_c$ denotes the polynomials defining the cross-ratio relations, say $z_{ik}(T)(1-z_{ij}(T))-1$ and $z_{il}(T)z_{ij}(T)-z_{ij}(T)+1$;
\item[---] $P_0 = \prod_{z \in \cZ} z(1-z)$ is the polynomial defining the forbidden values for the $z_{ij}(T)$.
\end{itemize}
We have to compute an exhaustive and exact description of the constructible set
$$
\cS = \{ z \in \cZ \ | \ P(z)=0, P \in L_e \cup L_f \cup L_h \cup L_c, \ P_0(z) \not = 0\}.
$$

Some preliminary remarks are essential in order to understand the contents of the present section.

For fixed holonomy values (for example unipotent solutions), a
straightforward approach that would consist in first solving the large
algebraic system given by the consistency equations, the holonomy
equations and the cross-ratio relation would certainly fail because of
the number of variables ($12 n$) and the degree. Moreover it may
happen that this system is not 0-dimensional even if the holonomies
are fixed.



We chose to perform exclusively exact computation, using elimination techniques.
We will first simplify the system using specific adapted pre-processing, then we will use general methods essentially based on Gröbner basis computations (see \cite{CLO1997}).

When the system we obtain is 0-dimensional, we then make use of
the Rational Univariate Representation (\cite{Rouillier1999}) to get
formal parameterizations of the solutions.



\subsection{The pre-processing}\label{sec-prepro}

Using the cross-ratio relations at each vertex of a tetrahedron one can reduce the number of variables to $4n$.
We pick a variable at each vertex and write the equations in terms of these variables.
For example, such a substitution
reduces the study of the $4_1$ knot complement to the solutions of a 0-dimensional system depending on $8$ variables
(instead of $24$).

After these simplifications, some new simple affine relations  $A z_{ij}+B$ may appear,
where $A$ divides some power of $P_0$,  $A$ being a product of $z$ or $(1-z)$ for some $z \in \cZ - \{ z_{ij} \}$. The coordinate $z_{ij}$ may be replaced everywhere by $-{B}/{A}$. We then obtain a simplified system.

This simplification pre-processing reduces the study of the $4_1$ knot complement to a simplified system of two equations in two variables. The other coordinates of the solutions can then be recovered from simple relations.

Note that the degree of the equations in the final simplified system depends strongly on the way these simplifications are performed.
Blind simplifications may lead to consider systems of equations with very large degrees. The general methods for solving such systems formally have polynomial complexity in the Bézout's bound. Unappropriate choices in the simplifications will generate huge systems that are impossible to solve with {\em state of the art} algorithms. 

We have implemented some tricky choice functions that minimize the degree of the final simplified system. There is no interest in describing them in detail here, but this straightforward step must be carefully implemented if one wants to succeed in the determination of the solutions.

\subsection{Solving the reduced system using Gröbner Basis}
We now have to solve a system of equations and (simple) inequalities.
In the case of systems with fixed holonomy values (for example
unipotent solutions), three cases may occur:
\begin{itemize}
\item[---]
The set of equations defines a 0-dimensional variety (finite set of points);
\item[---]
The set of equations defines a variety of positive dimension but the
full constructible set is 0-dimensional (consider for example the
system $\{X(XY-1)=0,X(X^2+Y^2-1)=0\}$ and the constructible set
$\{X(XY-1)=0,X(X^2+Y^2-1)=0,X\neq 0\}$).
\item[---]
The constructible set has positive dimension.
\end{itemize}
In the case of unipotent solutions for problems with 2 or 3
tetrahedra, the first case happened only for the {\tt m004} variety
($4_1$ complement), the last case occurs for example for the {\tt m003} variety ($4_1$-sister) and also for the {\tt m006}.
The second case occurs for example for the {\tt m007} and the {\tt m015} ($5_2$ complement).
The procedure we follow is divided into the following steps:
\begin{enumerate}
\item {\bf Compute a Gr\"obner basis} of the reduced system of equations. 
\item {\bf Saturate the ideal} by the polynomials defining the inequalities. 
\item {\bf Compute the dimension} of this ideal. 
\end{enumerate}
These computations are based on Gröbner Basis computations.

\def\cI{\cal I}
A Gröbner basis of a polynomial ideal ${\cI}$
is a set of generators of  ${\cI}$, such that
there is a natural and unique way (the {\em normal form}) of reducing canonically a polynomial $P \Mod{\cI}$.
A Gröbner basis is uniquely defined for a given admissible ordering on
the monomials.
Given
$\cI \subset \Q[Y_1,\ldots,Y_k][X_{1},\ldots,X_m]$,
the use of an {\em elimination ordering}
such that $Y_i<X_j$ for all $i,j$, allows to deduce
straightforwardly a Gr\"obner basis of $\cI \cap  \Q[Y_1,\ldots,Y_k]$ and therefore to eliminate $X_1, \ldots, X_m$.
(see \cite[Ch. 3]{CLO1997}).

We obtain the saturation of the ideal $\cI \subset \Q[X_1,\ldots,X_m]$ by a polynomial $P_0$ by  computing a Gröbner
basis of $\Bigl ( \cI+\langle TP_0-1\rangle \Bigr ) \cap  \Q[X_1,\ldots,X_m]$, where $T$ is a new
independent variable. 
This is the usual way for computing an ideal whose zeroes set is the Zariski closure of a given
constructible set (the constructible set itself when it is 0-dimensional).

Even if the basic principles used for computing a Gröbner basis are simple
(see \cite{CLO1997} for a general description),
their effective computation is known to be hard.
The first algorithm due to Buchberger has been a lot
improved and we use, in practice, the algorithm $F_4$ by J.-C. Faugère (see
{\cite{Faugere1999}}), implemented in recent versions of {\sc Maple}, and
known to be currently the fastest variant.

Once a Gr\"obner Basis of ${\cI}$ is known, one can compute
the associated Hilbert polynomial and then deduce the (Hilbert) dimension and
the (Hilbert) degree of ${\cI}$ (see \cite{CLO1997} Chapter 9 - Section
3).

\subsection{Prime decompositions}
In general,
the ideal generated by the equations as well as its saturation by the inequalities is
not prime. When it is of positive dimension, it has components of mixed dimensions.
A key point at this stage is to be able to compute a prime or at least a primary
decomposition.

There are several existing strategies for computing such a
decomposition and it would be too long to enumerate them in the
present article. One can mention two classes: methods based on
Gröbner bases with an elimination ordering and factorization (see
\cite{CLO1997} - chapter 7) and those based on triangular sets (see
\cite{ALM1999} for example). Note  however that the Maple function we use for
this operation implements heuristics in order to select a favorable strategy.

\subsection{Rational Univariate Representation}\label{sec:rur}

When a system is 0-dimensional, the so called Rational Univariate Representation (see
\cite{Rouillier1999}) defines a one-to-one correspondence between the solutions and the roots
of some univariate polynomial $f$, preserving the multiplicities.
If this systems depends on $m$ variables $X_1, \ldots ,X_m$,
the Rational Univariate Representation consists in a polynomial $f \in \Q[Y]$ and a set of
$m+1$ polynomials such that the solutions are
$$\{ x_i = \Frac{f_{i}(\gamma)}{f_0(\gamma)},\ i= 1, \ldots, m \ |\  f(\gamma)=0\}.$$

The algorithm is based on the fact that, when $\cI$ is 0-dimensional,
$V=\Q[X_1,\ldots,X_m]/\cI$ is a finite dimensional $\Q$-algebra 
generated by the monomials that are irreducible modulo $\cI$. 
These elements than can be directly obtained
from a Gröbner basis. 

There exists a linear combination $Y = \sum_{i=1}^m \lambda_i X_i \in V$ 
which {\em separates} the zeroes of $\cI$.   
Then $f$ is the characteristic polynomial of the map $V\rightarrow V$, $P \mapsto Y \cdot P$.

The polynomials $f, f_i, i = 0, \ldots, m$ are then computed by linear algebra.

This full parametrization is computed in two steps:
first compute a suitable set of generators for the ideal (the best being a Gröbner basis),
then compute a Rational Univariate Representation of the reduced
0-dimensional system.

Factorization of $f$ leads to consider that the solutions are parameterized by a finite number of algebraic numbers
$\gamma$, defined by their minimal polynomials $f$ and each coordinate $z$ is a polynomial $f_z (\gamma)$.

Then, we substitute the variables that have been extracted during the
pre-processing by their expression in $\Q(\gamma) = \Q[Y]/(f)$ using modular computation.
We obtain a rational  parametrization of the full system.
Each of these parametrization may be described as
$$\mathrm{rur}:=\{z=f_{z}(\gamma),\ z \in \cZ \ |\  f(\gamma)=0\}.$$
Furthermore, $\gamma \in \R$ if an only if $z(\gamma) = f_z(\gamma) \in \R$ for all $z \in \cZ$.

\subsection{Sorting the solutions}
In this section we explain an exact procedure to sort solutions in
$\SL(3,\bR)$, $\PGL(2,\bC)$ and $\PU(2,1)$ for 0-dimensional components.
Thanks to the previous computations, the solutions of the initial system are
expressed by means of rational parameterizations (here $f$ is a prime polynomial)
$$\mathrm{rur}:=\{z=f_{z}(Y),\ z \in \cZ\ |\ f(Y)=0\}.$$
\begin{enumerate}
\item Real solutions are in correspondence with the real zeroes of $f$.   The number of real roots can be determined exactly
by Sturm algorithms and their approximate values are computed by certified algorithms (see )

\item The hyperbolic solutions are extracted from the global solutions as follows.
Equation \ref{eq:hyp} becomes
$${z_{12}} \equiv {z_{21}} \equiv {z_{34}} \equiv {z_{43}} \Mod{f}.$$
From the Rational Univariate Representation it is just a matter of testing if
some formal coordinates are equal in a field extension.

\item For the CR solutions, we consider the $3\times n$ equations \ref{eq:cr}:
$$z_{ij} z_{ji} = \overline{z_{kl}}\, \overline{z_{lk}}.$$
Let $\gamma = x + iy$ be one of the root of $f$.
We obtain a new 0-dimensional system in $\Q[x,y]$ consisting of the $2$ polynomial equations:
$\Re f(\gamma)=0, \ \Im f(\gamma)=  0$
and the $2 \times 3n$ equations:
$$\Re  \left [{z_{ij}(\gamma)}{z_{ji}(\gamma)} \right ] =\Re \left [{{z_{kl}(\gamma)}} \, {{z_{lk}(\gamma)}}\right ],
\Im  \left [{z_{ij}(\gamma)}{z_{ji}(\gamma)} \right ] = - \Im \left [{{z_{kl}(\gamma)}} \, {{z_{lk}(\gamma)}}\right ].$$
This last system may be solved by computing a Rational Univariate Representation.
There exists an integer $\lambda$ such that
the algebraic number $\eta = x+ \lambda y$ describes the solutions:
$$\mathrm{rur}_{\mathrm{CR}}:=\left \{x=\Frac{g_{x}(\eta)}{g_{1}(\eta)},y=\Frac{g_{y}(\eta)}{g_{1}(\eta)}\ | \ g(\eta)=0 \right \}.$$
Real solutions of the system are in a one-to-one correspondance with the real roots of $g$.
\end{enumerate}
\subsection{Computing the rank of the boundary holonomy representation}
We compute first the meridian and the longitude as upper triangular matrices,
 say
$$g_M =
\left(
\begin{matrix}
1 & a & c\\
0 & 1 & b \\
0 & 0 & 1
\end{matrix} \right ), \quad
g_L =
\left(
\begin{matrix}
1 & a' & c'\\
0 & 1 & b' \\
0 & 0 & 1
\end{matrix} \right )
.$$
We easily obtain that $g_M^n =
\left(
\begin{matrix}
1 & n\,a & n\, c + \frac 12 n(n-1)ab\\
0 & 1 & n\, b \\
0 & 0 & 1
\end{matrix} \right )$.
We thus deduce that $g_M^n = g_L^{n'}$ iff $\Frac{a}{a'} = \Frac{b}{b'} = \Frac{2c - ab}{2c'-a'b'} = \Frac{n'}{n} \in \Q$.
This can be decided by computing first the reduced row echelon form  of $
\left ( \begin{matrix}
a & a' \\
b & b '\\
2c-ab &2c'-a'b'
\end{matrix} \right )
$
in the ground field $\Q[Y]/(f(Y))$ where $f$ is the minimal irreducible polynomial parameterizing the solutions.
If the reduced row echelon form has rank one, we test if $a$ and $a'$ are $\Q$-dependant or not.

\subsection{Certified numerical approximations}

For simplicity, and in order to get a human readable output, we often express the results by means of numerical approximations instead of formal (large) expressions. We assert that all these approximations are certified (in practice, all the digits are correct but the last one). We make use of multi-precision floating point numbers coupled with interval arithmetic or well known results on error control for the evaluation of univariate polynomials  that allows to accurately evaluate the Rational Univariate Representations at the roots of a given univariate polynomial. A fully treated example in the case of systems with two variables  is described in \cite{bouzidi:hal-00802698}.

%% file: calculs0-III.tex
\section{An example: the figure-eight knot}\label{sec:m004}
\begin{figure}[!ht]
\begin{center}
\includegraphics[scale=0.50]{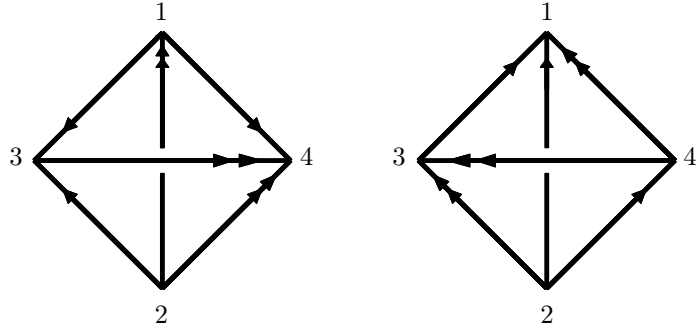}
\put(-215,105){1}
\put(-270,50){3}
\put(-160,50){4}
\put(-215,-10){2}
\put(-69,105){1}
\put(-125,50){3}
\put(-15,50){4}
\put(-69,-10){2}
\end{center}
\caption{Triangulation of {\tt m004}}
\label{fig:m004}
\end{figure}

The figure eight knot is given as a gluing of two tetrahedra as in Figure \ref{fig:m004}.  We refer to it
to set the equations of compatibility.

Let $z_{ij}$ and $w_{ij}$ be the coordinates associated to the edge $ij$ of each of the tetrahedra.


We obtain the edge and face equations:
$$
(L_e)
\left \lbrace
\begin{array}{rcl}
z_{{23}}w_{{24}}z_{{13}}w_{{31}}z_{{14}}w_{{21}}&=&1\\
z_{{32}}w_{{42}}z_{{31}}w_{{13}}z_{{41}}w_{{12}}&=&1\\
z_{{24}}w_{{23}}z_{{21}}w_{{43}}z_{{34}}w_{{41}}&=&1\\
z_{{42}}w_{{32}}z_{{12}}w_{{34}}z_{{43}}w_{{14}}&=&1.
\end{array}
\right.\quad
(L_f)
\left \lbrace
\begin{array}{rcl}
z_{{21}}z_{{31}}z_{{41}}w_{{21}}w_{{31}}w_{{41}}&=&1\\
z_{{12}}z_{{32}}z_{{42}}w_{{13}}w_{{23}}w_{{43}}&=&1\\
z_{{13}}z_{{23}}z_{{43}}w_{{12}}w_{{32}}w_{{42}}&=&1\\
z_{{14}}z_{{24}}z_{{34}}w_{{14}}w_{{24}}w_{{34}}&=&1.
\end{array}
\right.
$$
\subsection{Boundary holonomy representation}
\begin{figure}[ht]
\begin{center}
\begin{asy}
  // 
  // 
  size(12cm);
  defaultpen(1);
  usepackage("amssymb");
  import geometry;
  //arrowbar Fleche = Arrow(SimpleHead, 3, 10);
  // 
  point o = (0,0);
  pair ph = (2,0);
  pair pv = (1,2);

  // 
  draw(o -- o+pv -- o+ph -- (3,2) -- (4,0)--(5,2)-- o+3*ph --(7,2)--(8,0)--(9,2)--(1,2)--cycle, blue);
  draw(o -- o+4*ph, blue);

  // 
  // dot(o,black+5);
  // 
  label("{\small $z_{34}$}", o, 1.5*dir(35));
  label("\small $ z_{41}$", o+ph, 1.5*dir(35));
  label("\small $z_{12}$", o+2*ph, 1.5*dir(35));
  label("\small $z_{23}$", o+3*ph, 1.5*dir(35));

  label("\small $ z_{32}$", o+pv, 4*S);
  label("\small $ z_{42}$", o+pv+ph, 4*S);
  label("\small $z_{14}$", o+pv+2*ph, 4*S);
  label("\small $ z_{24}$", o+pv+3*ph, 4*S);

  label("\small $z_{31}$", o+ph, 1.5*dir(145));
  label("\small $z_{43}$", o+2*ph, 1.5*dir(145));
  label("\small $z_{13}$", o+3*ph, 1.5*dir(145));
  label("\small $ z_{21}$", o+4*ph, 1.5*dir(145));

  label("\small$w_{12}$", o+pv, 1.5*dir(-35));
  label("\small $w_{32}$", o+pv+ph, 1.5*dir(-35));
  label("\small $ w_{21}$", o+pv+2*ph, 1.5*dir(-35));
  label("\small $w_{41}$", o+pv+3*ph, 1.5*dir(-35));

  label("\small$ w_{13}$", o+ph, 4.5*N);
  label("\small $ w_{34}$", o+ph+ph, 4.5*N);
  label("\small $w_{24}$", o+ph+2*ph, 4.5*N);
  label("\small $ w_{43}$", o+ph+3*ph, 4.5*N);

  label("\small $ w_{14}$", o+pv+ph, 1.5*dir(-145));
  label("\small $ w_{31}$", o+pv+2*ph, 1.5*dir(-145));
  label("\small $w_{23}$", o+pv+3*ph, 1.5*dir(-145));
  label("\small $ w_{42}$", o+pv+4*ph, 1.5*dir(-145));

 // 
 // 
 //
  draw ((8,2)--(7.5,1)--(1.5,1)--(1,0),red,Arrow(5bp,position=1.4));

// 
draw ((.5,1)..(1,0),green, Arrow(5bp));
draw ((8,2)..(8.5,1),green,Arrow(5bp));

  // 
  // pair p1 = (4, 2);
  // pair p2 = (7, -1);
  // guide H = p1{E}..{S}p2;
  // draw(H,red);

\end{asy}
\caption{The holonomy of the figure eight cusp. The red line corresponds to $L$ and the green line to $M$} \label{fig:eight}
\end{center}
\end{figure}

The meridian and the longitude are computed from the Snappea triangulation.  We follow the algorithm to compute the holonomy along
the meridian and the longitude along appropriate paths:
\begin{itemize}
\item[]
$(421)^1\rightarrow\ (423)^1\rightarrow\ (324)^0\rightarrow\ (314)^0$
\item[]
$(421)^1\rightarrow\ (431)^1\rightarrow\ (214)^0\rightarrow\ (234)^0\rightarrow\ (243)^1\rightarrow\ (241)^1\rightarrow\ (134)^0
\rightarrow\ (132)^0\rightarrow\ (312)^1\rightarrow\ (342)^1\rightarrow\ (432)^0\rightarrow\ (412)^0\rightarrow\ (134)^1
\rightarrow\ (132)^1\rightarrow\ (312)^0\rightarrow\ (314)^0$
\end{itemize}
The corresponding matrices are:
$$g_M =
\left[ \begin {array}{ccc} {\Frac {z_{{12}}z_{{21}}}{z_{{34}}w_{{
24}}}}&-{\Frac {z_{{12}}z_{{13}}+z_{{34}}z_{{24}}}{z_{{34}}z_{{2
4}}z_{{13}}}}&-{\Frac {z_{{32}}w_{{42}}}{z_{{13}}}}
\\ \noalign{\medskip}0&1&z_{{32}}w_{{42}}\\ \noalign{\medskip}0&0&{
\Frac {w_{{42}}}{z_{{34}}}}\end {array} \right]
g_L = \left[
\begin {array}{ccc}
{\Frac {z_{{13}}z_{{31}}w_{{14}}w_{{23}}}{w_{{
31}}z_{{42}}w_{{42}}z_{{24}}}}&*&*\\
0&1&*\\ 0&0&{\Frac {z_{{31}}w_{{13}}z_{{13}}w_{{24}}}
{z_{{42}}w_{{32}}z_{{24}}w_{{41}}}}\end {array} \right]
$$
We therefore deduce the unipotent holonomy equations: $A=A^*=B=B^*=1$ where
$$
A=\Frac {z_{{12}}z_{{21}}}{z_{{34}}w_{{24}}},\
A^*=\Frac {w_{{42}}}{z_{{34}}},\
B={\Frac {z_{{13}}z_{{31}}w_{{14}}w_{{23}}}{w_{{31}}z_{{42}}w_{{42}}z_{{24}}}}, \
B^*={\Frac {z_{{31}}w_{{13}}z_{{13}}w_{{24}}}{z_{{42}}w_{{32}}z_{{24}}w_{{41}}}}.
$$

\subsection{Solutions}
Here the pre-processing provides a reduced system of 7 polynomial equations in the variables
$\{z_{14} , z_{43}\}$. All of the 22 other variables belong to $\Q(z_{14} , z_{43})$.
This reduced system is 0-dimensional and one can easily compute Rational Univariate Representations of its
zeroes.

One obtain four sets given by their minimal polynomials.
\begin{itemize}
\item[$\blacktriangleright$] $f_1 = Y^2-Y+1$. \\
The complete hyperbolic structure on the complement of the figure-eight knot is obtained from two conjugate solutions.
In fact, in that case, if $\omega^{\pm}=\frac{1\pm i\sqrt{3}}{2}$ is one root of $f_1$, then
$$
z_{12}=z_{21}=z_{34}=z_{43}=w_{12}=w_{21}=w_{34}=w_{43}=\omega^{\pm}
$$
is a solution of the equations as obtained in \cite{Thurston}. Its volume is $2.029883212\cdots$.

\item[$\blacktriangleright$] $f_2 = Y^2-Y+1$.
$$
z_{12}=z_{34}=w_{34}=w_{43}=
\overline z_{21}=\overline z_{43}=\overline w_{12}=\overline w_{21}
= \omega^{\pm}.
$$
This solution corresponds to a discrete representation
of the fundamental group of the complement of  knot in
$\PU(2,1)$ with faithful boundary holonomy.  Moreover, its action on complex hyperbolic space has limit set the full boundary sphere (\cite{falbeleight}).

\item[$\blacktriangleright$] $f_3 = Y^2+Y+2$.
Let $\gamma^{\pm} = -\frac 12 \pm i \frac 12 \sqrt 7$ be a root of $f_3$.
The coordinates are
$$
z_{21} = w_{21} = \overline{z}_{43} = \overline{w}_{12} = \frac{5-i\sqrt{7}}{4},\
z_{12} = w_{2,1}  = \frac{3-i\sqrt{7}}{8},
w_{34} = \overline{w}_{43} = -\frac{1+i\sqrt 7}2.
$$
\item[$\blacktriangleright$] $f_4 = 2Y^2+Y+1$.
$1/\gamma^{\pm} = -\frac 14 \pm i \frac 14 \sqrt 7$ is a root of $f_4$. We obtain
$$
z_{12} = \overline{z}_{34} =  \frac{3-i\sqrt{7}}{2},\
w_{43} = \overline{w}_{34} = \frac{5+i\sqrt{7}}{8},\
z_{21} = w_{21} = \overline{z}_{43} = \overline{w}_{12} =\frac{-1+i\sqrt{7}}{4}.
$$
\end{itemize}
All solutions corresponding to $f_2, f_3, f_4$ were obtained in
\cite{falbeleight} and $f_3, f_4$ correspond to spherical CR structures with unipotent boundary holonomy of rank one (\cite{DF}).
%
%

\subsection{Holonomy representation}

%
In order to obtain the representation of the fundamental group of the complement of the figure eight knot we identify the face
$(324)^0$ of the tetrahedron $0$ to the face
$(423)^1$ of the  tetrahedron $1$.  The generators are explicitly written from the remaining three face
identifications following  paths which pass through the chosen face
$(324)^0= (423)^1$ inside the tetrahedra and transport the projective basis
defined by each oriented face.

Here are the paths generating the identifications on the polyhedron:
\begin{enumerate}
\item $(134)^0\rightarrow\ (341)^0\rightarrow\ (342)^0\rightarrow\ (423)^0\rightarrow\ (234)^0=(243)^1\rightarrow\ (241)^1$
\item $(214)^0\rightarrow\ (142)^0\rightarrow\ (421)^0\rightarrow\ (423)^0\rightarrow\ (234)^0=(243)^1\rightarrow\ (432)^1\rightarrow\ (431)^1$
\item $(123)^0\rightarrow\ (231)^0\rightarrow\ (234)^0=(243)^1\rightarrow\ (432)^1\rightarrow\ (324)^1\rightarrow\ (321)^1$
\end{enumerate}

In order to obtain paths generating the fundamental group with a base point at face $(234)^0$,
the paths have to be conjugated, respectively, by
\begin{enumerate}
\item $(234)^0\rightarrow\ (342)^0\rightarrow\ (341)^0\rightarrow\ (413)^0\rightarrow\ (134)^0$
\item $(234)^0\rightarrow\ (342)^0\rightarrow\ (423)^0\rightarrow\ (421)^0\rightarrow\ (214)^0$
\item $(234)^0\rightarrow\ (231)^0\rightarrow\ (312)^0\rightarrow\ (123)^0$
\end{enumerate}
The coordinate transformations corresponding to each path are obtained
by multiplying the coordinate change matrices.  The generators of the
group representation are obtained as inverse matrices of the coordinate change matrices.

Applying the computations above to the solution $f_2$ we obtain the generators
$$
g_1=\left[ \begin {array}{rrr}
1&0&0\\
-w^{\pm}&1&0\\
- w^{\mp}&{-2w^{\mp}}&1
\end {array} \right],\,
g_2=\left[ \begin {array}{rrr}
3 w^{\mp}&6 w^{\mp}&1\\
3 w^{\pm}&-1&0\\
1&0&0
\end {array} \right], \,
g_3 = \left[ \begin {array}{rrr}
1& 2w^{\mp}&- w^{\mp}\\
0&1&{w^{\mp}}\\
0&0&1\end {array} \right].
$$
The group $\langle g_1, g_2, g_3 \rangle$ preserves the hermitian form
$
\left[ \begin {array}{rrr}
0&0&-{w^{\pm}}\\
0&1&0\\
-w^{\mp}&0&0
\end {array} \right]
$.

For the  solution $f_3$ we obtain

\noindent
{\small
$
g_1=\left[ \begin {array}{ccc}
1&0&0\\
-\Frac{2}{\gamma^{\pm}}&1&0\\
-\Frac{\gamma^{\pm}}{2}&\Frac{2}{\gamma^{\pm}-1}&1
\end {array} \right]
$,
$g_2= \left[ \begin {array}{ccc}
\Frac{8}{\gamma^{\pm}+3}&\Frac{3}{\gamma^{\pm}+2}&-\Frac{4}{3\gamma^{\pm}+2}\\
-\Frac{3}{\gamma^{\pm}}&1&0\\
\Frac{2}{\gamma^{\pm}+2}&0&0
\end {array} \right]$,
$g_3=\left[ \begin {array}{ccc}
1&\Frac{2}{\gamma^{\pm}+2}&-\Frac{1}{\gamma^{\pm}}\\
0&1&\gamma^{\pm}\\
0&0&1
\end {array} \right].
$
The group $\langle g_1, g_2, g_3 \rangle$ preserves the hermitian form
$
\left[ \begin {array}{rrr}
0&0&{\gamma^{\mp}}\\
0&1&0\\
\gamma^{\pm}&0&0
\end {array} \right]
$.
}

Finally, for the solution $f_4$ we obtain

\noindent
{\small
$
g_1 = \left[ \begin {array}{ccc}
1&0&0\\
\noalign{\medskip}
\Frac{2}{\gamma^{\pm}-2}&1&0
\\ \noalign{\medskip}
-\Frac{2}{\gamma^{\pm}}&\frac 12 \gamma^{\pm} -1&1\end {array}
 \right]$,
$
g_2= \left[ \begin {array}{ccc}
-2&-\frac 12 \gamma^{\pm}-3&-1\\ \noalign{\medskip}
\Frac{8}{\gamma^{\pm}+6}&1&0\\ \noalign{\medskip}
-1&0&0
\end {array} \right]$,
$g_3=\left[ \begin {array}{ccc}
1&\Frac{4}{\gamma^{\pm}+3}&-\gamma^{\pm}\\ \noalign{\medskip}
0&1&-\Frac{2}{\gamma^{\pm}-2}
\\ \noalign{\medskip}0&0&1\end {array} \right].
$
}
The group $\langle g_1, g_2, g_3 \rangle$ preserves the hermitian form
$
\left[ \begin {array}{rrr}
0&0&-\frac 12\\
0&1&0\\
-\frac 12&0&0
\end {array} \right]
$.

Observe that the three generators satisfy the relations
$$
g_3g_2^{-1}g_1^{-1}g_2=1 \ \ \ {\mbox{and}}\ \ \ g_3^{-1}g_2g_1^{-1}g_3g_1=1
$$
which give a presentation of the fundamental group of the figure eight knot.
\subsection{Boundary holonomy}

For all solutions in $\PU(2,1)$ we always obtain  $g_L = g_1^{-1} $.
The meridian is respectively
$$
g_M =
\left[ \begin {array}{rrr}
1&0&0\\
0&1&0\\
\pm\sqrt {3}&0&1
\end {array} \right],\
g_M = {g_L}^{-3}, \ g_M = {g_L}^3 .
$$

\subsection{Some non-unipotent configurations of the figure-eight knot}\label{sec:m004ni}

In this section we give examples of non-unipotent representations in order to show how our method can be applied without modification
to a general holonomy condition. We will only deal with a few examples for the figure eight knot and reserve a complete study of that knot for
a later paper.

We keep the notations of section \ref{sec:m004}:
$$
A=\Frac {z_{{12}}z_{{21}}}{z_{{34}}w_{{24}}},\
A^*=\Frac {w_{{42}}}{z_{{34}}},\
B={\Frac {z_{{13}}z_{{31}}w_{{14}}w_{{23}}}{w_{{31}}z_{{42}}w_{{42}}z_{{24}}}}, \
B^*={\Frac {z_{{31}}w_{{13}}z_{{13}}w_{{24}}}{z_{{42}}w_{{32}}z_{{24}}w_{{41}}}}.
$$
We will first impose that $A$ and $A^*$ are primitive cubic roots of unity,
and we find a finite set of solutions.
In the second case we will impose that $A$ and $A^*$ are primitive sixth roots of unity.
We find here a finite number of isolated solutions as well a one 1-dimensional component.

We start with the case when $A$ and $A^*$ are primitive cubic roots of unity.
Here the pre-processing provides a reduced system of 5 polynomial equations in the variables
$\{z_{14}, z_{24}, z_{43}, w_{14}\}$. All of the 20 other variables belong to $\Q(z_{14}, z_{24}, z_{43}, w_{14})$.
This reduced system is 1-dimensional. After saturation the system becomes 0-dimensional. It has 6 components.


Let $\zeta$ be a primitive twelfth root of unity.

The first component corresponds to representations  in $\PGL(2,\Q(\zeta^2))$ and we get $A=1/A^*= \zeta^4$, $B=B^*=1$.
The volume is 0, we have $g_M^3= g_L = 1$.

The second and the third component give 8 solutions  in $\PGL(3,\Q(\zeta))$. We get
$A=1/A^*= \zeta^4$, $B=B^*=1$. Their volumes are equal to 0 and we have $g_M^3= g_L = 1$.

The fourth and the fifth components give 8 solutions in $\PGL(3,\Q(\alpha))$ where $\alpha = \pm \sqrt 5 \pm i \sqrt 3$.
We get $A=A^*=\zeta^4,\, B=B^*=1$. In this case we have $g_M^3 = g_L$ or $g_M^3 = g_L^{-1}$.

The last component gives 4 solutions in $\PGL(3,\bC)$ and 4 solutions in $\PU(2,1)$.
We get $A=A^*=\zeta^4$ and $B=B^*$ is a root of
${X}^{4}-175\,{X}^{3}-327\,{X}^{2}-175\,X+1$.
This last polynomial has 2 conjugated roots with moduli 1 and two
real roots.  Therefore, the boundary holonomy group has rank 2.
These solutions belong to $\PU(2,1)$ when $\abs{B}=1$ and the four others give a volume
$.973235\cdots$, when $B$ is real.

\medskip\par\noindent
We now look at solutions with $A$ and $A^*$ primitive sixth roots of unity.
We find here seven 0-dimensional components and one 1-dimensional component.

There are 4 hyperbolic solutions corresponding to $A=1/A^*=\zeta^2$ and $B=1/B^* = 7 \pm 4 \sqrt 3$. Their volume is
$1.2212874588\cdots$. The meridian satisfies $g_M^6=1$.

There are two components in $\PGL(3,\Q(\zeta))$, satisfying $A=1/A^*=\zeta^2$, $B=B^*=-1$. They all satisfy
$g_M^3 = g_L, \, g_L^2= 1$. None of these belong to $\PU(2,1)$.

There are three components in $\PGL(3,\Q(\beta))$ where $\beta$ is a root of $X^4+X^3-X^2-X+1$. We have
$A=1/A^*=\zeta^2$, $B=B^*=-1$. They all satisfy
$g_M^3 = g_L, \, g_L^2= 1$. None of these belong to $\PU(2,1)$.

There is a component giving 4 solutions in $\PGL(3,\bC)$ and 4 solutions in $\PU(2,1)$.
We get $A=A^*=\zeta^4$ and $B=B^*$ is a root of
$X^4-27X^3+25X^2-27X+1$.
This last polynomial has 2 conjugated roots with moduli 1 and two
real roots.  Therefore, the boundary holonomy group has rank 2.
These solutions belong to $\PU(2,1)$ when $\abs{B}=1$ and the four others give a volume
$1.730258\cdots$, when $B$ is real.

The 1-dimensional component corresponds to reducible representations.
We have: $A=1/A^*=\zeta^2$ and $B=1/B^*=-1$.
The volume of these configurations is 0.
For all these configurations, the face variables are all equal to $-1$.
We obtain $z_{ij}z_{ji} = \zeta^{\pm 2},\, w_{ij}w_{ji} = \zeta^{\mp 2}$ for all $i<j$.
We always obtain $g_M^3=g_L$ and $g_L^2=1$.

\section{Description of the unipotent solutions for the first hyperbolic manifolds}

As explained in the last section we solve the system consisting of compatibility equations and the unipotent holonomy conditions plus the inequalities which prevent
the $z$-coordinates to be 0 or 1.   The solutions will, in turn, give rise to representations in $\PSL(3,\bC)$.

We have computed an exhaustive list of solutions for the eleven first 3-manifolds and for the Whitehead link complement.
Their main properties are summarized in the Table \ref{DesSols}.
For these examples the dimension of solutions is at most 1 and there are always isolated solutions besides the hyperbolic one.
Five of them are 0-dimensional:
\hyperref[sec:m004]{\tt m004}, \hyperref[sec:m007]{\tt m007}, \hyperref[sec:m009]{\tt m009}, \hyperref[sec:m015]{\tt m015} and
the \hyperref[sec:wh]{Whitehead link} complement.
All others are 1-dimensional:
\hyperref[sec:m003]{\tt m003}, \hyperref[sec:m006]{\tt m006}, \hyperref[sec:m010]{\tt m010}, \hyperref[sec:m011]{\tt m011},
\hyperref[sec:m016]{\tt m016}, \hyperref[sec:m017]{\tt m017}, \hyperref[sec:m019]{\tt m019}.
We further analyse the solutions to decide if the corresponding representation is in $\PSL_3(\bR)$, $\PSL_2(\bC)$, $\PSL_2(\bR)$ or $\PU(2,1)$ and compute their volume.  Observe also that, by Mostow rigidity, only one solution corresponds to the complete hyperbolic structure.  This solution is identified as the only one having all positive imaginary parts of the z-coordinates.

In Table \ref{DesSols}, we indicate the {\em name} of the variety, the number {\em 1-D} of 1-dimensional prime components,
the {\em degrees} of the prime 0-dimensional components.  The isolated solutions studied are the ones not contained in the 1-dimensional components. We indicate the total number of solutions such  that the corresponding representations are in $\PSL_{3}(\bC)$, $\PSL_3(\bR)$, $\PSL_2(\bC)$, $\PSL_2(\bR)$ and $\PU(2,1)$. Observe that $\PSL_3(\bR)\cap \PSL_2(\bC)=\PSL_2(\bR)$  and $\PSL_3(\bR)\cap \PU(2,1)=\PSL_2(\bR)$ so $\PSL_2(\bR)$ solutions are being also redundantly counted in $\PSL_3(\bR)$, $\PSL_2(\bC)$ and $\PU(2,1)$. We also indicate the different
positive volumes we obtain, writing in boldface the volumes of representations in $\PSL(2,\bC)$. Keep in mind though that real solutions as well as solutions giving rise to $\PU(2,1)$ representations have zero volume.

For the 0-dimensional case, we briefly describe the solutions and in some cases the group representations. In particular, we enumerate in the examples all cases when the boundary holonomy has rank one.  More details can be obtained in the website,
\href{https://who.rocq.inria.fr/Fabrice.Rouillier/SGT/Home_page.html}{.../SGT}
including matrix representations of the fundamental groups. In the following description, for the sake of simplicity, we will say that a solution is CR (short for Cauchy-Riemann), when
the corresponding representation is conjugated to a representation in $\PU(2,1)$. For  each variety we singled out a $\PU(2,1)$ representation with boundary holonomy
of rank one.  We think that they are natural candidates for being holonomies of a spherical CR structure on the variety.  Again, the difficulty in obtaining a CR structure lies in the definition of an appropriate 2-skeleton
 (see \cite{falbeleight} for a discussion).

The 1-dimensional components in the examples are not completely described.   They are all degenerate, that is, they are either reducible representations or they have trivial unipotent boundary holonomy and their image is a cyclic group (we thank M. Deraux, A. Guilloux and C. Zickert for discussions about these components which helped us to understand their meaning).  The reducible families arise from our definition of tetrahedra.  In fact we don't impose that the lines of the quadruple of flags be
in general position and the solutions with all the lines of all flags passing through a fixed point gives the family of reducible representations. We don't enumerate the
$\PSL(2,\bC),\PU(2,1)$ or $\PSL(3,\bC)$ solutions in these components.  But an important observation is that  in all examples up to three tetrahedra the 1-dimensional components contain at most a finite number of $\PSL(2,\bC)$ solutions which turn out to be real so the volumes presented in the table contain indeed  all volumes of $\PSL(2,\bC)$ representations arising from the given triangulation.
On the other hand, the 0-dimensional solutions are never reducible nor have trivial boundary holonomy.

\newcolumntype{M}[1]{>{\raggedleft}m{#1}}
{\small
\begin{table}[!ht]
\renewcommand{\arraystretch}{1.2}
\setlength{\tabcolsep}{2pt}
\begin{tabular}{||r|r|r|r|r|r|r|r|M{42mm}||}
\hline
 & & \multicolumn{6}{c}{0-dimensional prime components} &  \\
 \hline
 & & \multicolumn{6}{c}{Number(s) of Solutions} &\\
 \hline
Name&1-D&Ext. Degrees&$\PGL_{3}(\bC)$&$\PSL_{3}(\bR)$&$\PSL_2(\bC)$&$\PSL_2(\bR)$&$\PU(2,1)$&$Volumes$ \tabularnewline
 \hline
m003&	2&$2, 2, 8, 8$&	20&	0&	2&	0&	2&	  	  0.648847 {\bf 2.029883}\tabularnewline
 \hline
m004&	0&$2, 2, 2, 2$&	8&	0&	2&	0&	6&	  {\bf 2.029883} \tabularnewline
 \hline
m006&	2&	$6, 6, 12, 28$&	43&	1&	3&	1&	15&	  	  0.707031  0.719829  0.971648
  1.284485 {\bf 2.568971}\tabularnewline
 \hline
 m007&	0&	$3, 6, 8, 8, 8$&	33&	1&	3&	1&	15&	   0.707031  0.822744  1.336688 {\bf 2.568971}	\tabularnewline
 \hline
 m009&	0&	$2, 4, 4, 4, 6,8$&	28&	2&	2&	0&	8&	  	  0.507471  0.791583  1.417971 {\bf 2.666745}\tabularnewline
 \hline
 m010&	2&	$2, 6, 6,  12, 12$&	38&	0&	2&	0&	4&	  	  0.251617  0.791583  0.809805
  0.982389  1.323430 {\bf 2.666745} \tabularnewline
 \hline
 m011&	1&	$3, 4, 16, 64$&	87&	5&	7&	3&	21&	  	  0.226838  0.251809  0.328272
  0.397457  0.452710  0.643302
  0.685598  0.700395  0.724553
  0.770297  0.879768  {\bf 0.942707} 0.988006
  1.099133  1.184650  1.846570  {\bf  2.781834}\tabularnewline
 \hline
 m015&	0&	$3, 4, 4, 6, 6$&	23&	3&	3&	1&	11&	  	  0.794323  1.583167 {\bf 2.828122}\tabularnewline
 \hline
 m016&	1&	$3, 3, 10, 50$&	66&	4&	6&	4&	24&	  	  0.296355  0.403707  0.710033
  0.753403  0.773505  0.796590
  0.886451  1.135560  1.422985
  1.505989 {\bf 2.828122}\tabularnewline
 \hline
 m017&	3&	$3, 4, 6, 6, 44$&	63&	1&	3&	1&	21&	  	  0.527032  0.794323  0.801984
  0.828705  1.252969  1.588647 {\bf 2.828122}\tabularnewline
 \hline
 m019&	1&	$4, 4, 22, 84$&	114&	6&	8&	4&	24&	  	  0.027351  0.062112  0.323395
  0.332856  0.347159  0.411244
  0.467624  0.524801  0.544151
  0.599455  0.638404  0.738805
  0.758111  0.798098  0.851139
  0.916588  1.101800  1.130263
  1.190919  {\bf 1.263709}  1.340255  2.111776
  {\bf  2.944106}
  \tabularnewline
\hline
Wh. link & 0& $2, 2, 4, 4, 10,10$ &32 & 0& 2& 0& 14 &1.132196  1.683102 {\bf 3.663862}
  \tabularnewline
\hline
 \end{tabular}
 \vspace{.2in}
 \caption{Description of the solutions}
\label{DesSols}
\end{table}
 }
\newpage

\subsection{The variety {\tt m007} }\label{sec:m007}
There are three simplices with parameters $u_{ij}$, $v_{ij}$ and $w_{ij}$, $1\leq i,j\leq 4$.
\begin{figure}[!ht]
\begin{center}
\includegraphics[scale=0.50]{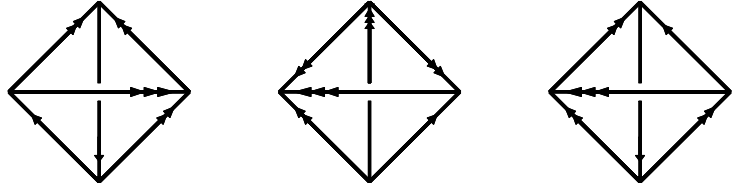}
\end{center}
\caption{Triangulation of {\tt m007}}
\label{fig:m007}
\end{figure}

\begin{itemize}
\item[$\blacktriangleright$] $f_1 = Y^3-2Y^2-1$. \\
The solutions giving rise to $\PSL_2(\bC)$ representations (containing the two conjugate complete hyperboilc ones) are
$$
\begin{array}{rcl}
u_{12} = u_{21} =u_{34} =u_{43} &=&\gamma,\\
v_{12} = v_{21} =v_{34} =v_{43} &=&2 \gamma - \gamma^2 = -1 /\gamma,\\
w_{12} = w_{21} =w_{34} =w_{43} &=&\gamma.
\end{array}
$$
where $\gamma$ is a root of $f_1 = Y^3 - 2Y^2-1$. The volume of the hyperbolic solution is
$\hbox{vol} = 2.5689706009\cdots$.
\item[$\blacktriangleright$] $f_2 =Y^6-Y^5+Y^4-2Y^3+Y^2+1$. \\
2 solutions are CR. The four others  define representations with volume
$.70703052208\cdots$.
\item[$\blacktriangleright$] $f_3 =Y^8-2Y^7+2Y^6-6Y^5+7Y^4-2Y^3+3Y^2-Y+1$. \\
4 solutions are CR. The volume of the other representations is
$0.82274406556\cdots$
\item[$\blacktriangleright$] $f_4 =3Y^8-6Y^7+7Y^6-4Y^5+4Y^4+4Y^3+4Y^2+1$. \\
4 solutions are CR. Their volume is
$.82274406556\cdots$
\item[$\blacktriangleright$] $f_5 =16Y^8-20Y^7+23Y^6-27Y^5+10Y^4-5Y^3+9Y^2+2Y+4$.\\
4 solutions are CR. The four others have volume
$1.3366875264\cdots$
\end{itemize}
In conclusion, there are 33 solutions,  15 being CR.

\subsection{The variety {\tt m009}}\label{sec:m009}
There are three simplices with parameters $u_{ij}$, $v_{ij}$ and $w_{ij}$, $1\leq i,j\leq 4$.
\begin{figure}[ht]
\begin{center}
\includegraphics[scale=0.50]{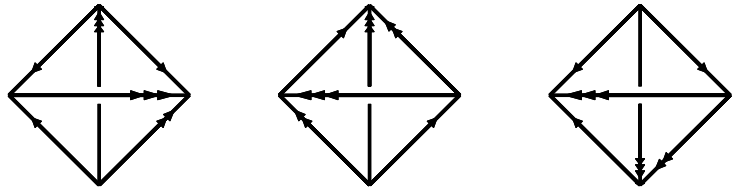}
\end{center}
\caption{Triangulation of {\tt m009}}
\label{fig:m009}
\end{figure}

We obtain
\begin{itemize}
\item[$\blacktriangleright$] $f_1 = Y^2-Y+2$. \\
The $\PSL(2,\bC)$ solutions (containing the two conjugate complete hyperbolic ones) are given by
$$
\begin{array}{rcl}
u_{12} = u_{21} =u_{34} =u_{43} &=&\gamma,\\
v_{12} = v_{21} =v_{34} =v_{43} &=&\frac 14 (\gamma+1)\\
w_{12} = w_{21} =w_{34} =w_{43} &=&\gamma.
\end{array}
$$
where $\gamma= \frac 12 (1\pm i \sqrt 7)$ is a root of $f_1$. The volume of the hyperbolic solution is
$\hbox{vol} = 2.6667447834\cdots$.
\item[$\blacktriangleright$] $f_2 = f_3 =  Y^4+Y^3-Y^2-Y+1$. \\
They correspond to 8 solutions with the same volume $0.50747080320\cdots$ .
\item[$\blacktriangleright$] $f_4 =  Y^4+2Y^3-Y^2-2Y-4 = (Y^2+Y-1 - \sqrt 5)(Y^2+Y-1+\sqrt 5)$. \\
There are two real solutions $\gamma^{\pm} = -\frac 12 \pm \frac 12 \sqrt{5+4\sqrt 5}$
and two conjugate CR solutions corresponding to $\gamma^{\pm} = -\frac 12 \pm \frac 12 i
\sqrt{-5+4\sqrt 5}.$
We have
$$
\begin{array}{l}
u_{12}=w_{34} = \Frac{\gamma^{\pm}+3}{\gamma^{\pm}+1},\,
u_{21}=w_{43} = \gamma^{\pm},\,
u_{34}=w_{12} =\Frac{\gamma^{\pm}-2}{\gamma^{\pm}},\\
u_{43}=w_{21} = - 1 - \gamma^{\pm},\,
v_{12}=v_{34} = \Frac{1}{\gamma^{\pm}+3},\,
v_{21}=v_{43} = \Frac{1}{2-\gamma^{\pm}}.
\end{array}
$$
The two real solutions give rise to representations in $\PSL(3\bR)\setminus \PSL(2,\bR)$.
Moreover, for all of these solutions we find that the meridian $g_M$ and the longitude $g_L$ satisfy $g_M g_L^2= 1$.
We find
$$
g_M =  \left[ \begin {array}{ccc}
1 & -2\,{\Frac {\gamma^{\pm}+2}{\gamma^{\pm}\, \left(\gamma^{\pm} +1\right) }}&
-2{\Frac {\gamma^{\pm}+2}{\gamma^{\pm}\, \left( \gamma^{\pm}+3 \right) }}\\\noalign{\medskip}
0 & 1 &\Frac{2}{3+\gamma^{\pm}}\\\noalign{\medskip}
0&0&1
\end {array}
\right].
$$
\item[$\blacktriangleright$] $f_5 =  Y^6-Y^4+2Y^3+Y^2-Y+1$. \\
There are two CR solutions and 4 solutions with the same volume $0.79158333031\cdots$.
\item[$\blacktriangleright$] $f_6 =  8Y^8-16Y^7+22Y^6-25Y^5+16Y^4-6Y^3+Y^2+3Y+1$. \\
There are 4 CR solutions and 4 solutions with the same volume $1.4179708859\cdots$.
\end{itemize}
In conclusion we found 28 solutions, 8 being CR.
There are 3 different volumes apart from the hyperbolic one.
\subsection{The $5_2$ knot complement}\label{sec:m015}
There are three simplices with parameters $u_{ij}$, $v_{ij}$ and $w_{ij}$, $1\leq i,j\leq 4$.
\begin{figure}[ht]
\begin{center}
\includegraphics[scale=0.50]{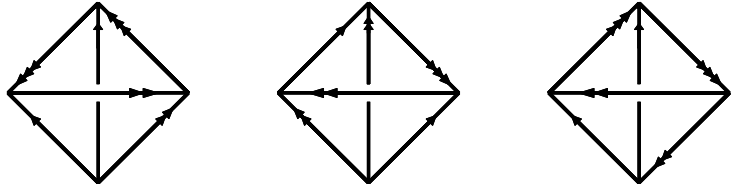}
\end{center}
\caption{Triangulation of {\tt m015}}
\label{fig:m015}
\end{figure}
We find 23 solutions in 5 sets of Galois conjugate solutions.
\begin{itemize}
\item[$\blacktriangleright$] $f_1 = Y^3-Y+1$.\\
The solutions correspond to the two conjugate complete hyperbolic solutions and a real one.
$$\begin{array}{rcl}
u_{12} = u_{21} =u_{34} =u_{43} &=& \gamma,\\
v_{12} = v_{21} =v_{34} =v_{43} &=& \gamma,\\
w_{12} = w_{21} =w_{34} =w_{43} &=& \gamma.
\end{array}
$$
where $\gamma$ is a root of $f_1$.  The volume is $2.828122\cdots$.
\item[$\blacktriangleright$] $f_2 =Y^4-Y^2+1$.\\
We obtain 4 solutions (the 4 12th-roots of unity) that are all CR.
%
\item[$\blacktriangleright$] $f_3 =4Y^4-4Y^3-Y^2+Y-1$.\\
We obtain 2 real solutions (giving representations in $\PSL(3,\bR)\setminus \PSL(2,\bR)$) that are the roots of
$$
4Y^2 -2Y - (1-\sqrt 5),$$
and 2  complex conjugate solutions that are the roots of
$$
4Y^2 -2 Y  - (1+\sqrt 5).
$$
These two last solutions are CR.
We obtain, for these four solutions:
$$
\begin{array}{r}
u_{21} = v_{21}=v_{34}=w_{21} = \Frac{1}{1-\gamma}\,\\
v_{12}=u_{43}=v_{43}=w_{43}=\Frac{2}{2\gamma +1},\\
u_{12}=w_{12}=-\gamma(\gamma-1),\,
u_{34}=w_{34}=\frac 14 - \gamma^2.
\end{array}
$$
Note that in this case the meridian and the longitude are equal:
$$
g_M = g_L= \left[ \begin {array}{ccc}
1 & \Frac{2}{(2\gamma-1)^2}&\Frac{1}{\gamma(2\gamma-1)^2}
\\\noalign{\medskip}
0 & 1 &\Frac{1}{\gamma}\\\noalign{\medskip}
0&0&1
\end {array}
\right].
$$
The holonomy group, in this case, is generated  by\\
\centerline{$g_1 = g_M^{-1}$,
$
g_2 = \left[ \begin {array}{ccc}
1&\Frac{2}{(2\gamma-1)^2}&\Frac{1}{\gamma(2\gamma-1)^2}\\\noalign{\medskip}
\Frac{2}{1-2\gamma} &\Frac{1-4\gamma}{1-4\gamma^2}&\gamma-1 \\\noalign{\medskip}
\Frac{1}{2\gamma^2(1-2\gamma)}&\Frac{\gamma-1}{2\gamma}&\gamma
\end {array}
\right],
$}
\centerline{$
g_3 = \left[ \begin {array}{ccc}
1&0&0\\\noalign{\medskip}
-\Frac{1}{2\gamma-1} & 1&0 \\\noalign{\medskip}
\Frac{1}{2\gamma^2(1-2\gamma)}& -\Frac{1}{2\gamma^2}&1
\end {array}
\right],
$
$
g_4=
\left[ \begin {array}{ccc}
\Frac{1+4\gamma}{2\gamma}&\Frac{1}{2\gamma(2\gamma-1)}&\Frac{2\gamma}{1-2\gamma}\\\noalign{\medskip}
\Frac{1}{2\gamma} &\Frac{2\gamma-1}{2\gamma}&0\\\noalign{\medskip}
1&0&0\end {array}
\right].
$}
where $\gamma$ is one of the four roots of $f_3$.
\item[$\blacktriangleright$] $f_4 =Y^6-Y^5-3Y^3+2Y^2+Y+1$.\\
Only two of these roots correspond to CR structures.
The four others have same volume: $1.583167\cdots$.
\item[$\blacktriangleright$] $f_5 =17Y^6-5Y^5+33Y^4-11Y^3+26Y^2-4Y+8$.\\
Only two of these roots correspond to CR structures. The others define representations of volume
$0.794323\cdots$.
\end{itemize}
In conclusion we obtain 23 solutions and 11 of them are CR.
\subsection{The Whitehead link}\label{sec:wh}
\begin{figure}[ht]
\begin{center}
\includegraphics[scale=0.50]{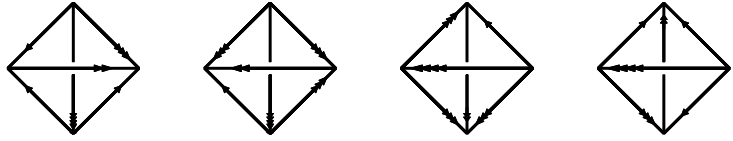}
\end{center}
\caption{Triangulation of Whitehead link complement}
\label{fig:wh}
\end{figure}

Here we have four tetrahedra with parameters $u_{ij}$, $v_{ij}$, $w_{ij}$ and $x_{ij}$ with $1\leq i,j\leq 4$ glued as in Figure \ref{fig:wh}.

There are 6 different groups of solutions.
\begin{itemize}
\item[$\blacktriangleright$] $f_1=Y^2+1$.
The conjugate complete hyperbolic solutions are given by
$$
u_{12}=u_{21}=u_{34}=u_{43}=v_{12}=v_{21}=v_{34}=v_{43} =\pm i,
$$
$$
w_{12}=w_{21}=w_{34}=w_{43}=w_{12}=w_{21}=w_{34}=w_{43} =\pm i.
$$
\item[$\blacktriangleright$] $f_2=Y^2+Y+4$.\\
We obtain
$$
\begin{array}{l}
u_{12}=v_{34}=w_{12}=x_{34}=\frac 38 - \frac 18 \gamma,\\
u_{21}=v_{43}=w_{21}=x_{43}=\gamma,\\
u_{34}=v_{12}=w_{34}=x_{12}=\overline{\gamma},\\
u_{43}=v_{21}=w_{43}=x_{21}=\frac 38 - \frac 18 \overline{\gamma}.
\end{array}
$$
where $\gamma = -\frac 12 \pm \frac 12 i \sqrt{15}$.

The holonomies of the meridian and longitude on the first torus (it corresponds to the vertex 1 of the first tetrahedron) are both equal to
$$
\left[ \begin {array}{ccc}
1&2&1\\\noalign{\medskip}
0&1&1\\\noalign{\medskip}
0&0&1\end {array}
\right],
$$
while the holonomies of the meridian and longitude of the second torus are equal to
$$
 \left[ \begin {array}{ccc} \frac 14+\frac 14\,i\sqrt {15}&\frac 34+\frac 34\,i\sqrt {15}&
\frac 32+\frac 12\,i\sqrt {15}\\ \noalign{\medskip}
-\frac 38-\frac 18\,i\sqrt {15}&-2-\frac 12
\,i\sqrt {15}&-{\frac {21}{8}}-\frac 38\,i\sqrt {15}\\ \noalign{\medskip}
\frac 32 &{\frac {21}{4}}+\frac 14\,i\sqrt {15}&{\frac {19}{4}}+\frac 14\,i\sqrt {15}
\end {array} \right].
$$
The representation corresponding to this solution was also obtained in \cite{PW} where it is shown that it is the holonomy of a spherical CR
structure on the complement of the Whitehead link.
\item[$\blacktriangleright$] $f_3=f_4=Y^4+Y^3+2Y^2-Y+1$.
All of these 8 solutions are CR.
\item[$\blacktriangleright$] $f_5=f_6=Y^{10}-4Y^8+Y^7+7Y^6-4Y^5-8Y^4+4Y^3+5Y^2+1$.
4 of these 20 solutions are CR-spherical.
\end{itemize}
Among these 32 solutions, 14 are CR.